\documentclass[12pt]{article}
\usepackage{amsfonts}
\usepackage{amssymb}
\usepackage{amsthm}
\usepackage{amsmath}
\usepackage{latexsym}
\usepackage{mathrsfs}
\usepackage{color}
\usepackage{bm}
\begin{document}
\newtheorem{Lemma}{\quad Lemma}
\newtheorem{pron}{\quad Proposition}
\newtheorem{thm}{\quad Theorem}
\newtheorem{Corol}{\quad Corollary}
\newtheorem{exam}{\quad Example}
\newtheorem{defin}{\quad Definition}
\newtheorem{remark}{\quad Remark}
\newcommand{\la}{\frac{1}{\lambda}}
\newcommand{\sectemul}{\arabic{section}}
\renewcommand{\theequation}{\sectemul.\arabic{equation}}
\renewcommand{\thepron}{\sectemul.\arabic{pron}}
\renewcommand{\theLemma}{\sectemul.\arabic{Lemma}}
\renewcommand{\thethm}{\sectemul.\arabic{thm}}
\renewcommand{\theCorol}{\sectemul.\arabic{Corol}}
\renewcommand{\theexam}{\sectemul.\arabic{exam}}
\renewcommand{\thedefin}{\sectemul.\arabic{defin}}
\renewcommand{\theremark}{\sectemul.\arabic{remark}}

\title{Asymptotics of convolution with the semi-regular-variation tail and its application to risk
\thanks {Research supported by National Science Foundation of China
(No.~10671139,11601043) } }
\author{\small Zhaolei Cui\ \ Edward Omey\ \ Wenyuan Wang\ \ Yuebao Wang
\thanks{Corresponding author. Telephone: 86 512
67422726. Fax: 86 512 65112637. E-mail: ybwang@suda.edu.cn}\\
{\small\it 1). School of Mathematics and Statistics, Changshu Institute of Technology,}\\
{\small\it Suzhou, P. R. China, 215500} \\
{\small\it 2). Faculty of Economics and Business-Campus Brussels, KU Leuven,}\\
{\small\it Warmoesberg Brussels, Belgium, 26, 1000} \\
{\small\it 3). School of Mathematics, Xiamen University, Xiamen, P. R. China, 361005}\\
{\small\it 4). School of Mathematics, Soochow University, Suzhou, P. R. China, 215006}\\
}

\date{}

\maketitle {\noindent\small {\bf Abstract }}\\

{\small In this paper, according to a certain criterion, we divide the exponential distribution class into \textcolor[rgb]{0.00,0.07,1.00}{some} subclasses. One of them is closely related to the regular-variation-tailed distribution class, so it is called the semi-regular-variation-tailed distribution class. In the \textcolor[rgb]{1.00,0.00,0.00}{new} class, although all distributions are not convolution equivalent,
\textcolor[rgb]{1.00,0.00,0.00}{there} still have some good properties. We give the precise tail asymptotic expression of \textcolor[rgb]{1.00,0.00,0.00}{convolutions} of these distributions, and prove that \textcolor[rgb]{1.00,0.00,0.00}{the class} is closed under convolution. In addition, we do not need to require the corresponding random variables to be identically distributed. Finally, we apply these results to a discrete time risk model with stochastic returns, \textcolor[rgb]{1.00,0.00,0.00}{thus} obtain the precise asymptotic estimation of the finite time ruin probability.\\ 



\section{\bf Introduction}

The theory and applications of the class of subexponential distributions has been
studied intensively, see, for example, \textcolor[rgb]{1.00,0.00,0.00}{Chistyakov (1964),} Embrechts et al. (1997) and Foss et al. (2013).
These distributions have a number of attractive properties. Among
others they obey the principle of a \textcolor[rgb]{0.00,0.07,1.00}{single} jump\textcolor[rgb]{1.00,0.00,0.00}{, that is} 
the tail of the distribution of the sum of independent and identically distributed
(i.i.d.) random variables (r.v.'s) is \textcolor[rgb]{1.00,0.00,0.00}{asymptotically} equal to the tail of the largest value of the summands. The distribution of the sum of independent r.v.'s is called convolution, and the distribution of their product is called product convolution. In probability theory, the tailed analysis of the convolution and product convolution is an important topic, \textcolor[rgb]{1.00,0.00,0.00}{which} have important applications in various fields such as finance, insurance, queuing systems and so on. Since the product of positive r.v.'s can be transformed to a sum of
r.v.'s (by taking logarithms), the convolution of r.v.'s and subexponential distribution, as well as more general convolution equivalent \textcolor[rgb]{1.00,0.00,0.00}{distributions} are especially important.

In practice, many distributions are not convolution equivalent. Even in the common exponential distribution class, \textcolor[rgb]{1.00,0.00,0.00}{which contains the convolution equivalent distribution class,} there are many such distributions, we refer to Pitman (1979), Embrechts and Goldie (1980), Murphree (1989), Leslie (1989), Lin and Wang (2012), Wang et al. (2016), among others. In the present paper, we also provide some examples. 

In this way, some interesting questions arise naturally:

\textbf{Problem A.} Is there a precise asymptotic expression for the tail of convolution of distributions except for the convolution equivalent distributions?

\textbf{Problem B.} Can we remove the restriction that these corresponding r.v.'s are identically distributed?

\textbf{Problem C.} What are the relationship and the difference between the such distributions and other distributions?

In the conclusion of the paper, we should reply on these questions.
We restrict our objects of study to the exponential distribution class. In this class, except for the convolution equivalent distribution, we expect to find some distributions with such good properties. It is well known that the tail of a distribution belonging to the class is asymptotically equivalent to a function $e^{-\alpha x}f(x)$ for some non-negative constant $\alpha$, where the function $f\circ\ln$ is slowly varying\textcolor[rgb]{0.00,0.07,1.00}{, for example, see (1.1) of Kl\"{u}ppelberg (1989)}. In
\textcolor[rgb]{1.00,0.00,0.00}{this} paper we consider \textcolor[rgb]{0.00,0.07,1.00}{some} different assumptions about $f$.
According to the criterion, we divide the 
class into \textcolor[rgb]{0.00,0.07,1.00}{two mutually disjoint} subclasses, see Subsection 1.1 below. We are mainly concerned with one
of these subclasses. Although these distributions in the \textcolor[rgb]{1.00,0.00,0.00}{new} subclass are not convolution equivalent,
they still have other good properties. We obtain precise tail asymptotic
expressions for convolution of these distributions and prove that the subclass
is closed under convolution. In addition, we do not require that the underlying r.v.'s are identically distributed.

In the following, we introduce the research objects and the main results of this paper, respectively.

\subsection{\bf Concepts of some functions and distributions}

Firstly, we give some \textcolor[rgb]{1.00,0.00,0.00}{notations} and conventions around the regular-variation function.
Without \textcolor[rgb]{1.00,0.00,0.00}{further comment, all limits are} as $x\to\infty$. Let
$u$ and $v$ be two eventually positive functions. We write $u(x)\sim
v(x)$, if $\lim u(x)/v(x)=1$; we write $u(x)=o\big(v(x)\big)$, if $\lim
u(x)/v(x)=0$; we write $u(x)=O\big(v(x)\big)$, if $\limsup u(x)/v(x)<\infty$; we write
$u(x)\asymp v(x)$, if $u(x)=O\big(v(x)\big)$ and $v(x)=O\big(u(x)\big)$.

A measurable function $u$ supported on $\mathbb{R}$ is said to be regularly varying at $\infty$ with \textcolor[rgb]{1.00,0.00,0.00}{parameter} $\alpha\in\mathbb{R}$, denoted by $u\in\mathcal{R}_{\alpha}$, if $u(x)$ is \textcolor[rgb]{1.00,0.00,0.00}{eventually positive} and for all $t\in\mathbb{R^+}$,
\begin{eqnarray*}\label{regular-def-1}
u(xt)\sim t^\alpha u(x).
\end{eqnarray*}
If $\alpha=0$, the function is called slowly varying. 

We say that a function $u$ belongs to the exponential function class $\mathcal{L}_\alpha$
for some $\alpha\in\mathbb{R^+}\cup\{0\}$, if the compound function $u\circ\ln\in\mathcal{R}_{-\alpha}$, or equivalently, if for all $t\in\mathbb{R}$,
$$u(x-t)\sim e^{\alpha t}u(x).$$
In particular, \textcolor[rgb]{0.00,0.07,1.00}{the class $\mathcal{L}_0$ is called the long-tailed function class}.

According to the uniform convergence theorem, see Theorem 1.5.2 of Bingham et al. (1987), we know that, if the function $u\in\mathcal{L}_\alpha$, then
\begin{eqnarray*}
\mathcal{H}(u,\alpha)&=&\{h:\ \mathbb{R^+}\cup\{0\}\longmapsto\mathbb{R^+},h(x)\uparrow\infty,\textcolor[rgb]{0.00,0.07,1.00}{h(x)/x}\rightarrow0\ \text{and}\nonumber\\
&&\ \ \ \ \ \ \ \  \ \ \ u(x-t)\sim e^{\alpha t}u(x)\ \text{uniformly for all}\ \mid t\mid\le h(x)\}\neq\phi.
\end{eqnarray*}
The property is used \textcolor[rgb]{0.00,0.07,1.00}{in different situations, see for example Corollary 2.5 of Cline and Samorodnitsky (1994) and
Proposition 2 and Proposition 6 of Asmussen et al. (2003)}.

We say that a function $u$ belongs to the function class $\mathcal{S}_\alpha$
for some $\alpha\in\mathbb{R^+}\cup\{0\}$, if $u\in\mathcal{L}_{\alpha}$, $m(u,\alpha)=\int_0^\infty e^{\alpha y}u(y)dy<\infty$ and for all $x\in\mathbb{R}^+\cup\{0\}$,
$$u^{\otimes2}(x)=u\otimes u(x)\sim 2m(u,\alpha)u(x),$$
where $$u\otimes v(x)=\int_0^x u(x-y)v(y)dy$$
for two positive functions $u$ and $v$ on $\mathbb{R}$. In particular, the class $\mathcal{S}_0$ is called \textcolor[rgb]{1.00,0.00,0.00}{the} subexponential function class.

In this paper, let $U$ be a proper distribution supported on $\mathbb{R}$. If its tail $\overline{U}=1-U\in\mathcal{R}_{-\alpha}$
for some $\alpha\in\mathbb{R^+}\cup\{0\}$, then we say that the distribution belongs to the regular-variation-tailed distribution class,
denoted by $U\in\mathcal{R}(\alpha)$.

We say that a distribution $U$ belongs to the exponential distribution
class $\mathcal{L}(\alpha)$ for some $\alpha\in\mathbb{R^+}\cup\{0\}$, if $\overline{U}\in\mathcal{L}_\alpha$.
Particularly, the class $\mathcal{L}(0)$ is called \textcolor[rgb]{1.00,0.00,0.00}{the} long-tailed distribution class,
denoted by $\mathcal{L}$. 
Clearly, if $U\in\mathcal{L}(\alpha)$, then $\mathcal{H}(\overline{U},\alpha)\neq\phi$.

We say that a distribution $U$ belongs to the convolution equivalent distribution
class $\mathcal{S}(\alpha)$ for some $\alpha\in\mathbb{R^+}\cup\{0\}$, if $U\in\mathcal{L}(\alpha)$, $M(U,\alpha)=\int_{-\infty}^\infty e^{\alpha y}U(dy)<\infty$ and
$$\overline{U^{*2}}(x)\sim2M(U,\alpha)\overline{U}(x),$$
\textcolor[rgb]{1.00,0.00,0.00}{where $U*V$ is the convolution of two distributions $U$ and $V$, and $U^{*2}=U*U$.}
Similarly, we denote the class $\mathcal{S}(0)$ by $\mathcal{S}$, which is called the subexponential distribution class.


Now, we try to find a criterion for the classification of the class $\mathcal{L}(\alpha)$ for some $\alpha\in\mathbb{R}^+$.

Let $X$ be a r.v. with the distribution $V$ belonging to the class $\mathcal{L}(\alpha)$ for some $\alpha\in\mathbb{R^+}\cup\{0\}$ and  let $Y=e^X$ with a distribution $F$ on $\mathbb{R}^+$. Then we know that $F\in\mathcal{R}(\alpha)$, so that there exist a function $l\in\mathcal{R}_0$ on $\mathbb{R}^+$ such that, for all $x\in\mathbb{R}$,
\begin{eqnarray}\label{1000}
\overline{V}(x)=P(Y>e^x)=\overline{F}(e^x)\sim e^{-\alpha x}l(e^x)=e^{-\alpha x}f(x).
\end{eqnarray}
Clearly, 
\textcolor[rgb]{0.00,0.07,1.00}{the function $f$ (=$l\circ\exp$) on $\mathbb{R}$ or equivalently, the function $l$, is locally bounded,} i.e. for any $x_0\in\mathbb{R}^+$, there is a constant $C=C(f,x_0)>0$ such that $f(x)\le C$ for all $0\le x\le x_0$.

The class $\mathcal{L}(\alpha)$ can be divided into \textcolor[rgb]{0.00,0.07,1.00}{two disjoint subclasses as follows:
\begin{eqnarray*}
\mathcal{L}_1(\alpha)=
\Big\{V\in\mathcal{L}(\alpha): \int_0^\infty f(y)dy=\infty\ \text{in}\ (\ref{1000})\Big\}
\end{eqnarray*}
and
\begin{eqnarray*}
\mathcal{L}_2(\alpha)=\mathcal{L}(\alpha)\setminus\mathcal{L}_1(\alpha)=\Big\{V\in\mathcal{L}(\alpha):\ \int_0^\infty f(y)dy<\infty\ \text{in}\ (\ref{1000})\Big\}.
\end{eqnarray*}
In this paper, we focus on a subset of the class $\mathcal{L}_{1}(\alpha)$ that
\begin{eqnarray*}
\mathcal{L}_{11}(\alpha)=\cup_{\gamma\ge-1}\mathcal{L}_{11}(\alpha, \gamma)=\cup_{\gamma\ge-1}\Big\{V\in\mathcal{L}_1(\alpha):\  f\in\mathcal{R}_\gamma\ \text{in}\ (\ref{1000})\Big\}.
\end{eqnarray*}}

From the definition of the subclass $\mathcal{L}_{11}(\alpha)$,  we can see the close relation between the subclass and the regular-variation-tailed function class \textcolor[rgb]{1.00,0.00,0.00}{$\mathcal{R}_\gamma$}. So, we might as well call the former \textcolor[rgb]{1.00,0.00,0.00}{the} semi-regular-variation-tailed distribution class.

In Subsection 4.1, we give some specific examples to show that the \textcolor[rgb]{1.00,0.00,0.00}{classes $\mathcal{L}_{11}(\alpha)$ and $\mathcal{L}_{1}(\alpha)\setminus\mathcal{L}_{11}(\alpha)$} contain many distributions with natural shapes.

If a distribution $V\in\textcolor[rgb]{1.00,0.00,0.00}{\mathcal{L}_1(\alpha)}$ for some $\alpha\in\mathbb{R}^+$, then $\int_0^\infty f(y)dy=\infty$ if and only if $M(V,\alpha)=\infty$. In fact, for any two constants $0<c<s$ large enough, using integration by parts, we have
\begin{eqnarray*}
\textcolor[rgb]{0.00,0.07,1.00}{\int_0^s}e^{\alpha y}U(dy)&=&-e^{\alpha s}\overline{U}(s)+\overline{U}(0)+\alpha\Big(\int_0^c+\int_c^s\Big)e^{\alpha y}\overline{U}(y)dy\\
&\le&\overline{U}(0)+\alpha\int_0^ce^{\alpha y}\overline{U}(y)dy+2\alpha\int_c^sf(y)dy.\\
\end{eqnarray*}
Thus, \textcolor[rgb]{1.00,0.00,0.00}{if $M(V,\alpha)=\infty$, then $\int_0^\infty f(y)dy=\infty$}.
Conversely, by $f\in\textcolor[rgb]{0.00,0.07,1.00}{\mathcal{L}_0}$ and (\ref{thm05}) below, 
we have
\begin{eqnarray*}
\int_c^se^{\alpha y}U(dy)&\ge&-2f(s)+\alpha\int_c^sf(y)dy\Big/2
\ge\alpha\int_c^sf(y)dy\Big/4.
\end{eqnarray*}
Therefore, $M(V,\alpha)=\infty$ follows from $\int_0^\infty f(y)dy=\infty$.

In addition, if $\gamma>-1$, then $\int_0^\infty f(y)dy=\infty$ \textcolor[rgb]{1.00,0.00,0.00}{holds automatically}. On the other hand, if $\gamma<-1$, then $\int_0^\infty f(y)dy<\infty$.

Correspondingly, for \textcolor[rgb]{1.00,0.00,0.00}{a parameter} $\alpha\in\mathbb{R}^+$, the class $\mathcal{R}(\alpha)$ can also be divided into \textcolor[rgb]{0.00,0.07,1.00}{two} disjoint subclasses.
To this end, let $Y$ be a r.v. with the distribution $F_0\in\mathcal{R}(\alpha)$. We denote its positive part by
$Y^+=Y\textbf{1}(Y>0)$ with the distribution $F$ on $\mathbb{R}^+$, then for all $x\in\mathbb{R}$,
\begin{eqnarray}\label{10}
\overline{F}(x)=\overline{F_0}(x)\textbf{1}(x\ge0)+\textbf{1}(x<0),
\end{eqnarray}
thus $F\in\mathcal{R}(\alpha)$. Further, let r.v. $X=\ln Y^+$, then its distribution $V=F\circ \exp$ belongs to the class $\mathcal{L}(\alpha)$. For $\textcolor[rgb]{0.00,0.07,1.00}{1\le i\le2}$, we \textcolor[rgb]{1.00,0.00,0.00}{define the following classes:}
$$\mathcal{R}_i(\alpha)=\{F:\ V=F\circ \exp\in\mathcal{L}_i(\alpha)\}.$$
\textcolor[rgb]{0.00,0.07,1.00}{The two} classes do not intersect each other,
$$\mathcal{R}(\alpha)=\textcolor[rgb]{0.00,0.07,1.00}{\mathcal{R}_1(\alpha)\cup\mathcal{R}_2(\alpha)}\subset\mathcal{S}
\subset\mathcal{L},$$
and $V\in\mathcal{L}_i(\alpha)$ if and only if  $F\in\mathcal{R}_i(\alpha)$ \textcolor[rgb]{0.00,0.07,1.00}{for $i=1,2$. Among them,
\begin{eqnarray*}
\textcolor[rgb]{0.00,0.07,1.00}{\mathcal{R}_1(\alpha)}\supset\{F_0\in\mathcal{R}(\alpha):\ V=F\circ\exp\in\mathcal{L}_{11}(\alpha)\}=\mathcal{R}_{11}(\alpha)
\end{eqnarray*}}
and
\begin{eqnarray*}
\textcolor[rgb]{0.00,0.07,1.00}{\mathcal{R}_2(\alpha)}\supset\{F_0\in\mathcal{R}(\alpha):\ V=F\circ\exp\in\mathcal{S}(\alpha)\}=\mathcal{R}^*(\alpha),
\end{eqnarray*}
which is called strongly regular-variation-tailed distribution class introduced by Definition 2.1 of Li and Tang (2015).

In the above-mentioned definitions and properties of the class $\mathcal{L}(\alpha)$, when
$\alpha\in\mathbb{R}^+$ and the distribution $U$ is lattice, all variables
and constants should be restricted to values of the lattice span in the distribution $U$ and the function $f$, see Bertoin and Doney (1996).
For example, let $U$ be a lattice distribution on $\mathbb{N}$. We say that $U\in\mathcal{L}(\alpha)$, if for all $i\in\mathbb{N}$,
\begin{eqnarray*}
\overline{U}(k-i)\sim e^{\alpha i}\overline{U}(k)\ \ \big(thus\ f(k-i)\sim f(k)\big)\ \ as\ k\to\ \infty.
\end{eqnarray*}
Further, $U\in\mathcal{L}_{11}(\alpha)$ if and only if $U\in\mathcal{L}(\alpha),\ \sum_{i=0}^\infty f(i)=\infty$ and
\begin{eqnarray*}
f(ki)\sim i^\alpha f(k)\ \ as\ k\to\ \infty.
\end{eqnarray*}
In the above sense, we still call the function $f$ belongs to the class \textcolor[rgb]{0.00,0.07,1.00}{$\mathcal{L}_0$} or the class $\mathcal{R}_\alpha$.

Now we return to the class $\mathcal{L}(\alpha)$ for some $\alpha\in\mathbb{R}^+$.
\textcolor[rgb]{0.00,0.07,1.00}{Since the introduction of exponential distribution class by Chover et al. (1973a,b), the tail asymptotics of convolution and the closure under convolution or convolution roots for convolution equivalents distribution have been thoroughly studied, see, for example, Embrechts and Goldie (1982), Kl\"{u}ppelberg (1989), Pakes (2004), Foss and Korshunov (2007) and Watanabe (2008). However, there are few corresponding results on other exponential distributions.}
\textcolor[rgb]{1.00,0.00,0.00}{We are mainly concerned with the class $\mathcal{L}_1(\alpha)$ and its subclass $\mathcal{L}_{11}(\alpha)$. We have not found a study of the former.}
In the research related to \textcolor[rgb]{0.00,0.07,1.00}{the latter, Hashova and Li (2013) give the tail asymptotic expression of convolution of some special distributions in the subclass. And then}, Hashova and Li (2014) obtained an asymptotic expression for the finite time ruin probability in a discrete time risk model with both insurance risks and financial risks, natural logarithms of which also follow some special distributions in the \textcolor[rgb]{0.00,0.07,1.00}{subclass}, see Theorem A and Remark \ref{remark30} below. On the other hand,
Omey et al. (2017) systematically studied properties of some r.v.'s, which \textcolor[rgb]{1.00,0.00,0.00}{are called semi-heavy tailed}, from the two angles of the probability density for an absolutely continuous distribution and the probability sequence for a lattice distribution, respectively. For example, the paper considered such probability density $w(x)=e^{-\alpha x}f(x)$ for some positive constant $\alpha$ and regular-variation function $f$. If $f$ is a tail distribution, then $w$ also is a special semi-regular-variation-tailed distribution, because the function $f$ must be non-increasing and tend to zero.

Inspired by the \textcolor[rgb]{0.00,0.07,1.00}{three papers above, in a different way, we try to study the tail} asymptotic properties of convolution with a unified form for the distributions in the class $\mathcal{L}_1(\alpha)$ and its subclass $\mathcal{L}_{11}(\alpha)$.


\subsection{\bf Main results}

We first give a general result for the class $\mathcal{L}_1(\alpha)$ \textcolor[rgb]{1.00,0.00,0.00}{with some parameter} $\alpha\in\mathbb{R}^+$. Then, based on this, we give another \textcolor[rgb]{1.00,0.00,0.00}{precise} result for the class $\mathcal{L}_{11}(\alpha)$. \textcolor[rgb]{0.00,0.07,1.00}{
We assume that all lattice distributions in this paper are supported on $\mathbb{N}$.}

\begin{thm}\label{thm1}
For some integer $n\ge1$ and all $1\le i\le n+1$, let $Y_i$ be a r.v. with distribution $V_i\in\mathcal{L}_1(\alpha)$ for some $\alpha\in\mathbb{R}^+$. Further, assume that $Y_i,1\le i\le n+1$, are independent of each other. Then \textcolor[rgb]{1.00,0.00,0.00}{$V_1*\cdots *V_n\in\mathcal{L}_1(\alpha)$ and}
\begin{eqnarray}\label{170618-corol1-1}
\overline{V_1*\cdots *V_n}(x)\sim a^{n-1} e^{-\alpha x}f_1\otimes \cdots\otimes f_n(x)=o\big(\overline{V_1*\cdots *V_{n+1}}(x)\big),
\end{eqnarray}
where $a=\alpha$, when the distributions $V_i,1\le i\le n+1$ are non-lattice; $a=e^\alpha-1$, when they are lattice; $a=(e^\alpha-1)^{m/(n-1)}\alpha^{(n-m-1)/(n-1)}$, when $n\ge2$ and there is an integer $1\le m\le n-1$ such that the distribution $V_i$ is lattice for $1\le i\le m$ and the distribution $V_j$ is non-lattice for $m+1\le j\le n$. And for any $1\le j\le n+1$,
\begin{eqnarray}\label{170618-corol1-100}
P\big(\max\{Y_i:1\le i\le n+1\}>x\big)
=o\Big(\sum_{1\le i\neq j\le n+1}\overline{V_i*V_j}(x)\Big).
\end{eqnarray}
\end{thm}

\begin{remark}\label{remark10}
i) 
It is easy to see that \textcolor[rgb]{0.00,0.07,1.00}{for $\alpha>0$, the classes $\mathcal{L}_1(\alpha)$ and  $\mathcal{L}_{11}(\alpha)$ are} closed under convolution, that is if $V_1$ and $V_2$ belong to the class, then $V_1\ast V_2$ still belongs to the same one and has a \textcolor[rgb]{0.00,0.07,1.00}{heavier} tail than $V_i$ for $i=1,2$.
However, it is well known that, the class $\mathcal{S}(\alpha)$ does not have the property. This is an essential distinction between the two.

ii) For some $n\ge2$ and all $1\le i\le n$, we know that, if $V_i=V\in\mathcal{S}(\alpha)$ with some $\alpha\in\mathbb{R}^+\cup\{0\}$, then
\begin{eqnarray}\label{170618-thm1-101}
\overline{V^{*n}}(x)=P\Big(\sum_{i=1}^nY_i>x\Big)\sim M(V,\alpha)P\big(\max\{Y_i:1\le i\le n\}>x\big)\sim nM(V,\alpha)\overline{V}(x).
\end{eqnarray}
In the case that $\alpha=0$, these distributions or the corresponding r.v.'s obey the principle of ``a single big jump" or ``Max-Sum equivalence". And the second result of the theorem shows (\ref{170618-thm1-101}) does not hold for the distribution $V\in\mathcal{L}_1(\alpha)$ for some $\alpha\in\mathbb{R}^+$. This is another essential difference. 
However, there still exists another precise equivalent relation for convolution of distributions in the class $\mathcal{L}_1(\alpha)$, as the first expression in (\ref{170618-corol1-1}).

\end{remark}


In the following, \textcolor[rgb]{0.00,0.07,1.00}{denotes $\Gamma(\cdot)$ the Gamma-function such that $\Gamma(x)=\int_0^\infty y^{x-1}e^{-y}dy$ and $B(\cdot,\cdot)$ the Beta-function such that $B(\gamma_1+1,\gamma_2+1)=\int_0^1(1-y)^{\gamma_1}y^{\gamma_2}dy$}. Let $f$ be a non-negative function on $\mathbb{R}$ such that $f(x)>0$ eventually. For all $x\in\mathbb{R}^+\cup\{0\}$, we write
$$f^I(x)=\int_0^xf(y)dy.$$
We know that, if $f\in\mathcal{L}_0$, then by Proposition 1.5.9a of Bingham et al. (1987),
\begin{eqnarray}\label{thm05}
f(x)=o\big(f^I(x)\big).
\end{eqnarray}
\textcolor[rgb]{1.00,0.00,0.00}{Further, if $f\in\mathcal{R}_\gamma$ for some $\gamma\ge-1$, then by Proposition 1.5.9a and Theorem 1.5.11 of Bingham et al. (1987), $f^I\in\mathcal{R}_{\gamma+1}$ and}
\begin{eqnarray}\label{thm050}
\lim xf(x)/f^I(x)=\gamma+1.
\end{eqnarray}

\begin{thm}\label{corol1}
Under the conditions of Theorem \ref{thm1}, we further assume that the distribution $V_i\in\mathcal{L}_{11}(\alpha)$ for any $n\ge1$ and all $1\le i\le n+1$. \textcolor[rgb]{1.00,0.00,0.00}{Then $V_1*\cdots *V_{n+1}\in\mathcal{L}_{11}(\alpha)$ and the following results hold:}

i) If $f_i\in\mathcal{R}_{\gamma_i}$ and $\gamma_i>-1$, $1\le i\le n+1$, then
\begin{eqnarray}\label{170618-corol1-2}
\overline{V_1*\cdots *V_{n+1}}(x)&\sim& \prod_{j=1}^{n}B\Big(\sum_{k=1}^j\gamma_k+j,\gamma_{j+1}+1\Big)a^{n} e^{-\alpha x}x^{n}\prod_{i=1}^{n+1}f_i(x)\nonumber\\
&\sim&B\Big(\sum_{k=1}^n\gamma_k+n,\gamma_{n+1}+1\Big)a e^{\alpha x}x\overline{V_{n+1}}(x)\overline{V_1*\cdots *V_{n}}(x).
\end{eqnarray}

ii) If $f_i\in\mathcal{R}_{-1}$, $1\le i\le n+1$, then
\begin{eqnarray}\label{170709-corol1-1}
\overline{V_1*\cdots *V_{n+1}}(x)&\sim&a^{n} e^{-\alpha x}\sum_{i=1}^{n+1} f_i(x)\prod_{1\le j\neq i\le n+1}f_j^I(x).
\end{eqnarray}

iii) If there is some integer $1\le m\le n$ such that $f_i\in\mathcal{R}_{-1}$ for $1\le i\le m$, and $f_i\in\mathcal{R}_{\gamma_i}$ with some $\gamma_i>-1$ for $m+1\le i\le n+1$,  then
\begin{eqnarray}\label{170706-corol1-1}
&&\overline{V_1*\cdots *V_{n+1}}(x)\sim a^{n} e^{-\alpha x}x^{n-m}\nonumber\\
&&\ \ \ \ \ \ \ \ \cdot\prod_{i=m+1}^{n+1}f_i(x)\prod_{j=m+1}^{n}
B\Big(\sum_{s=m+1}^j\gamma_s+j-m,\gamma_{j+1}+1\Big)\prod_{r=1}^{m}f_r^I(x).
\end{eqnarray}
\end{thm}

\begin{remark}\label{remark100}
i) In the above two theorems, we only require that all distributions belong to the class $\mathcal{L}_1(\alpha)$ or $\mathcal{L}_{11}(\alpha)$ for some common parameter $\alpha\in\mathbb{R}^+$, because we do not need to consider other situations. In fact, for example, if $V_i\in\mathcal{L}_{11}(\alpha_i)$ with some $\alpha_i\in\mathbb{R}^+$ for $i=1,2$, where $\alpha_2>\alpha_1=\alpha$, then $\overline{V_2}(x)=o\big(\overline{V_1}(x)\big)$, $M(V_2,\beta)<\infty$ for each $\alpha\le\beta<\alpha_2$, and by Lemma 2.1 of Pakes (2004), we have
$$\overline{V_1*V_2}(x)\sim M(V_2,\alpha)\overline{V_1}(x).$$
Thus, $V_1*V_2\in\mathcal{L}_{11}(\alpha)$.

ii) In particular, if $V_i=V$ with $\gamma_i=\gamma$ and $f_i=f$ for all $1\le i\le n+1$, then we have,
for example, in (\ref{170709-corol1-1}),
$$\overline{V^{*(n+1)}}(x)\sim (n+1)a^{n} e^{-\alpha x}f(x)\big(f^I(x)\big)^{n}\sim (n+1)a^{n}\big(f^I(x)\big)^{n}\overline{V}(x).$$
For semi-heavy-tailed r.v.'s with common probability density function or lattice sequence, Omey et al. (2017) also obtained some corresponding results. 

iii) For the class $\mathcal{L}_1(\alpha)\setminus\mathcal{L}_{11}(\alpha)$, we also have some individual results, see Proposition \ref{pron41} and Proposition \ref{pron42} below.
\end{remark}

The rest of this paper consists of three sections. In Section 2, we give the proofs of Theorem \ref{thm1} and Theorem
\ref{corol1}. Further, we apply the obtained results to a
discrete time risk model with stochastic returns, \textcolor[rgb]{1.00,0.00,0.00}{thus} get the asymptotic
estimation of the finite time ruin probability in Section 3. Finally, in Section 4, we give some specific  
distributions with normal shapes and good properties in the classes $\mathcal{L}_{11}(\alpha)$ and $\mathcal{L}_1(\alpha)\setminus\mathcal{L}_{11}(\alpha)$ for some $\alpha\in\mathbb{R}^+$.
%
%

\section{\bf Proofs of the main results}
\setcounter{equation}{0}\setcounter{Lemma}{0}\setcounter{thm}{0}\setcounter{remark}{0}\setcounter{exam}{0}


\subsection{\bf Proof of Theorem \ref{thm1}}

In order to prove Theorem \ref{thm1},  we first give a lemma.

\begin{Lemma}\label{lem5}
i) For $i=1,2$, if $f_i\in\mathcal{L}_0$, then $f_1\otimes f_2\in\mathcal{L}_0$. Further, if $\int_0^\infty f_2(y)dy=\infty$, then \begin{eqnarray}\label{lemma2000}
f_1(x)=o\big(f_1\otimes f_2(x)\big).
\end{eqnarray}
ii) For $i=1,2$, let $X_i$ be a r.v. with distribution $V_i$ on $\mathbb{R}$, and they are independent of each other. Further, if $V_i\in\mathcal{L}_1(\alpha)$ for some $\alpha\in\mathbb{R}^+$, \textcolor[rgb]{1.00,0.00,0.00}{then
for any $c\in\mathbb{R}$},
\begin{eqnarray}\label{lemma201}
\overline{V_1*V_2}(x)&\sim&P(X_1+X_2>x,X_1>c,X_2>c)\sim P(X_1+X_2>x,X_i>c).
\end{eqnarray}
\end{Lemma}

The first result in i) is a minor generalization of Theorem 3 (b) of Embrechts and Goldie (1980) for $\alpha=0$ and Proposition 5 of Asmussen et al. (2003), where $f_i$ is a tail distribution for the former or a \textcolor[rgb]{0.00,0.07,1.00}{density} for the latter, $i=1,2$.

\proof i) Firstly, we prove the first conclusion. For any constant $t$, since $f_i\in\mathcal{L}_0$ for $i=1,2$, there is a function \textcolor[rgb]{0.00,0.07,1.00}{$h$ such that, $2h\in\mathcal{H}(f_1,0)\bigcap\mathcal{H}(f_2,0)$, and} when $x$ is large enough,
\begin{eqnarray*}
f_1\otimes f_2(x+t)&=&\Big(\int_0^{x-h(x)}+\int_{x-h(x)}^{x+t}\Big)f_1(x+t-y)f_2(y)dy=T_1(x)+T_2(x).
\end{eqnarray*}
Because \textcolor[rgb]{0.00,0.07,1.00}{$f_1\in\mathcal{L}_0$, both of $h$ and $h_1$ belong to $\mathcal{H}(f_1,0)\bigcap\mathcal{H}(f_2,0)$, where $h_1(x)=h(x)+t$,}
\begin{eqnarray*}
T_1(x)\sim \int_0^{x-h(x)}f_1(x-y)f_2(y)dy
\end{eqnarray*}
and
\begin{eqnarray*}
T_2(x)= \int_0^{h_1(x)}f_1(y)f_2(x+t-y)dy\sim\int_0^{h_1(x)}f_1(y)f_2(x-y)dy.
\end{eqnarray*}
And by (\ref{thm05}), we have
\begin{eqnarray*}
\int_{h(x)}^{h_1(x)}f_1(y)f_2(x-y)dy\sim f_2(x)f_1(h(x))t=o\big(T_2(x)\big).
\end{eqnarray*}
Therefore,
\begin{eqnarray*}
T_2(x)\sim \int_0^{h(x)}f_1(y)f_2(x-y)dy=\int_{x-h(x)}^x f_1(x-y)f_2(y)dy.
\end{eqnarray*}
So, $f_1\otimes f_2(x+t)\sim f_1\otimes f_2(x)$, that is $f_1\otimes f_2\in \mathcal{L}_0$.

Further, for any constant $x_0>0$, since $f_1\in\mathcal{L}_d$, we have
$$f_1\otimes f_2(x)\ge\int_0^{x_0}f_1(x-y)f_2(y)dy\sim f_1(x)\int_0^{x_0}f_2(y)dy.$$
\textcolor[rgb]{1.00,0.00,0.00}{Since $x_0$ is arbitrary and $\int_0^{\infty}f_2(y)dy=\infty$, relation (2.1) follows.}

ii) \textcolor[rgb]{1.00,0.00,0.00}{According to (\ref{170618-corol1-1}) in Theorem \ref{thm1}, we have
\begin{eqnarray}\label{lemma200}
\max\big\{\overline{V_1}(x),\overline{V_2}(x)\big\}=o\big(\overline{V_1*V_2}(x)\big).
\end{eqnarray}
For any $c\in\mathbb{R}$ and $1\le i\le2$}, because
\begin{eqnarray*}
P(X_1+X_2>x)-P(X_1+X_2>x,X_i\le c)\le P(X_i>x-c)\sim e^{\alpha c}\overline{V_i}(x),
\end{eqnarray*}
(\ref{lemma201}) follows from (\ref{lemma200}).
$\hfill\Box$\\

Now, we prove Theorem \ref{thm1}. \textcolor[rgb]{1.00,0.00,0.00}{By $\int_0^\infty f_1(y)dy=\infty$ and (\ref{lemma2000}), $\int_0^\infty f_1\otimes\cdots\otimes f_n(y)dy=\infty$, that is $V_1*\cdots *V_n\in\mathcal{L}_1(\alpha)$. To prove (\ref{170618-corol1-1})}, for a distribution $U$, we denote
$$\overline{U_I}(x)=\int_x^\infty\overline{U}(y)dy$$
for $x\in\mathbb{R}^+\cup\{0\}$. If $U\in\mathcal{L}(\alpha)$ for some $\alpha\in\mathbb{R}^+$, then
$$\overline{U_I}(x)=O\Big(\int_x^\infty e^{-\alpha y}l(e^y)dy\Big)=O\Big(\int_{e^x}^\infty z^{-\alpha-1}l(z)dz\Big)<\infty,$$
and by Karamata theorem, 
\begin{eqnarray}\label{thm04}
\overline{U_I}(x)\sim\overline{U}(x)/\alpha.
\end{eqnarray}

We first deal with the case $n=2$ and assume that $V_i$  is continuous for $i=1,2$.
Since $V_1,V_2$ and $V_1*V_2$  belong to $\mathcal{L}(\alpha)$, there exists a function $h\in\mathcal{H}(\overline{V_1},\alpha)\cap\mathcal{H}(\overline{V_2},\alpha)\cap\mathcal{H}(\overline{V_1*V_2},\alpha)$. \textcolor[rgb]{1.00,0.00,0.00}{Using this function $h$, (\ref{lemma201}) with $c=0$ and partial integration, for $x$ large enough we have}
\begin{eqnarray}\label{170616-thm-0}
&&\overline{V_1*V_2}(x)\sim P(X_1+X_2>x,X_1>0,X_2>0)\nonumber\\
&=&\Big(\int_{0}^{h(x)}+\int_{h(x)}^{x}\Big)\overline{V_1}(x-y)V_2(dy)+\overline{V_2}(x)\overline{V_1}(0)\nonumber\\
&=&\int_{0}^{h(x)}\overline{V_1}(x-y)V_2(dy)+\int_{0}^{x-h(x)}\overline{V_2}(x-y)V_1(dy)
+\overline{V_1}\big((x-h(x)\big)\overline{V_2}\big(h(x)\big)\nonumber\\
&=&I_1(x)+I_2(x)+I_3(x).
\end{eqnarray}

We first deal with $I_1(x)+I_3(x)$. Using (\ref{thm04}) for $U=V_1$ and integration by parts again, and note that $\overline{V_2}$ is a continuous function such that \textcolor[rgb]{0.00,0.07,1.00}{$\frac{d}{dy}(V_2)_I=\overline{V_2}(y)$} for all $y\in\mathbb{R^+}$, we have
\begin{eqnarray}\label{170616-thm-1}
&&I_1(x)+I_3(x)\sim\alpha\int_{0}^{h(x)}\overline{(V_1)_I}(x-y)V_2(dy)+\overline{V_1}\big((x-h(x)\big)\overline{V_2}\big(h(x)\big)\nonumber\\
&=&-\alpha \overline{(V_1)_I}(x-h(x))\overline{V_2}(h(x))+\alpha\overline{V_2}(0)\overline{(V_1)_I}(x)
+\alpha \int_{0}^{h(x)}\overline{V_1}(x-y)\overline{V_2}(y)dy\nonumber\\
&& +\overline{V_1}\big((x-h(x)\big)\overline{V_2}\big(h(x)\big).
\end{eqnarray}
Because $f_2$ is locally bounded, there exists a constant $C>0$ such that
$f_2(x)\le Ce^{\alpha x}\overline{V_2}(x)$ for all $x\in\mathbb{R}^+\cup\{0\}$. Further, by (\ref{thm04}) and (\ref{thm05}), we find that
\begin{eqnarray}\label{170616-thm-2}
\overline{(V_1)_I}\big(x-h(x)\big)\overline{V_2}\big(h(x)\big)
&\sim& \alpha^{-1}\overline{V_1}\big(x-h(x)\big)\overline{V_2}\big(h(x)\big)\ \ \big(=\alpha^{-1} I_3(x)\big)\nonumber\\
&\sim& \alpha^{-1}\overline{V_1}(x)\textcolor[rgb]{0.00,0.07,1.00}{e^{\alpha h(x)}e^{-\alpha h(x)}}f_2\big(h(x)\big)\nonumber\\
&=&o\Big(\overline{V_1}(x)\int_{0}^{h(x)}f_2(y)dy\Big)\nonumber\\
&=&o\Big(\int_{0}^{h(x)}\overline{V_1}(x-y)\overline{V_2}(y)dy\Big).
\end{eqnarray}
By (\ref{170616-thm-1}), (\ref{170616-thm-2}) and (\ref{thm04}),
\begin{eqnarray}\label{170616-thm-4}
I_1(x)+I_3(x)&\sim&\alpha \int_{0}^{h(x)}\overline{V_1}(x-y)\overline{V}_2(y)dy\textcolor[rgb]{1.00,0.00,0.00}{+\overline{V_2}(0)\overline{V_1}(x)}\nonumber\\
&\sim&\alpha e^{-\alpha x}\int_{0}^{h(x)}f_1(x-y)e^{\alpha y}\overline{V_2}(y)dy\textcolor[rgb]{1.00,0.00,0.00}{+\overline{V_2}(0)e^{-\alpha x}f_1(x)}\nonumber\\
&\sim&\alpha e^{-\alpha x}f_1(x)\int_{0}^{h(x)}e^{\alpha y}\overline{V}_2(y)dy.
\end{eqnarray}
For each fixed integer $n\ge1$, since $\int_{n}^{h(x)}e^{\alpha y}\overline{V_2}(y)dy\uparrow\infty$, there exists a constant $x_n\in\mathbb{R}^+$ such that
$$\max\Big\{n\int_{0}^{n}f_2(y)dy,n\int_{0}^{n}e^{\alpha y}\overline{V_2}(y)dy\Big\}
\le\min\Big\{\int_{n}^{h(x)}f_2(y)dy,\int_{n}^{h(x)}e^{\alpha y}\overline{V_2}(y)dy\Big\}$$
for all $x\ge x_n$. Without loss of generality, we set $nx_n\le x_{n+1}$ for all $n\ge1$. Let $h_1$ be a positive function on $\mathbb{R}^+\cup\{0\}$ such that
$$h_1(x)=\sum_{n=1}^\infty n^{1/2}\textbf{1}(x_n\le x<x_{n+1})+\textbf{1}(0\le x<x_1).$$
Then $h_1(x)\uparrow\infty$, $h_1(x)=o\big(h(x)\big)$ and
\begin{eqnarray}\label{170616-thm-41}
\max\Big\{\int_{0}^{h_1(x)}f_2(y)dy,\int_{0}^{h_1(x)}e^{\alpha y}\overline{V_2}(y)dy\Big\}
=o\Big(\min\Big\{\int_{n}^{h(x)}f_2(y)dy,\int_{n}^{h(x)}e^{\alpha y}\overline{V_2}(y)dy\Big\}\Big).
\end{eqnarray}
By (\ref{170616-thm-4}) and (\ref{170616-thm-41}), we have
\begin{eqnarray}\label{170616-thm-42}
I_1(x)+I_3(x)&\sim&I_1(x)\sim\alpha \overline{V_1}(x)\int_{h_1(x)}^{h(x)}e^{\alpha y}\overline{V_2}(y)dy\nonumber\\
&\sim&\alpha \overline{V_1}(x)\int_{h_1(x)}^{h(x)}f_2(y)dy\nonumber\\
&\sim&\alpha \overline{V_1}(x)\int_{0}^{h(x)}f_2(y)dy\nonumber\\
&\sim&\alpha\int_{0}^{h(x)}e^{-\alpha y}\overline{V_1}(x-y)f_2(y)dy\nonumber\\
&\sim&\alpha e^{-\alpha x}\int_{0}^{h(x)}f_1(x-y)f_2(y)dy.
\end{eqnarray}

In a similar way, we conclude that
\begin{eqnarray}\label{1001}
I_2(x)&\sim& \alpha\int_{h(x)}^x\overline{V}_1(x-y)\overline{V}_2(y)dy\sim \alpha e^{-\alpha x}\int_{h(x)}^xf_1(x-y)f_2(y)dy.
\end{eqnarray}

Hence, by (\ref{170616-thm-0}), (\ref{170616-thm-42}) and (\ref{1001}), it holds that
\begin{eqnarray}\label{170616-thm-5}
\overline{V_1*V_2}(x)\sim \alpha\int_0^x\overline{V_1}(x-y)\overline{V_2}(y)dy\sim\alpha e^{-\alpha x}f_1\otimes f_2(x).
\end{eqnarray}

Now, for $i=1,2$, we deal with the situation that $V_i$ doesn't have to be continuous. Let $Z_i$ be a random variable with a uniform distribution supported on $(0,\varepsilon)$ for some positive constant $\varepsilon$, where, when $V_i$ is lattice, $\varepsilon$ is the step size \textcolor[rgb]{1.00,0.00,0.00}{$1$}. We denote the distribution of $Y_i+Z_i$ by $G_i$. \textcolor[rgb]{1.00,0.00,0.00}{Assume} that $Y_i,Z_i,i=1,2$ are independent of each other.
Clearly, $G_i$ is absolutely continuous, and
\begin{eqnarray*}\label{40101}
G_i(x)&=&\int_0^\varepsilon V_i(x-y)dy/\varepsilon=\int_{x-\varepsilon}^x V_i(z)dz/\varepsilon=\int_{-\infty}^x\big(V_i(z)-V_i(z-\varepsilon)\big)dz/\varepsilon.
\end{eqnarray*}
It follows that $G_i$ has a density $g_i$ given by
$$g_i(x)=\big(V_i(x)-V_i(x-\varepsilon)\big)/\varepsilon=\big(\overline{V_i}(x-\varepsilon)-\overline{V_i}(x)\big)/\varepsilon\ \ a.s.$$

In the following, we first consider the situation that $V_i$ is non-lattice for $i=1,2$. By $V_i\in\mathcal{L}(\alpha)$ for some $\alpha>0$, we have
$$g_i(x)\sim\overline{V_i}(x)(e^{\alpha\varepsilon}-1)/\varepsilon=a(\varepsilon)\overline{V_i}(x).$$
Thus $g_i\in\mathcal{L}_\alpha$ and
\begin{eqnarray}\label{170616-thm-50}
\overline{G_i}(x)\sim a(\varepsilon)\overline{V_i}(x)/\alpha\sim a(\varepsilon)e^{-\alpha x}f_i(x)/\alpha=e^{-\alpha x}f_{0i}(x).
\end{eqnarray}
Further, as the proof of (\ref{170616-thm-42}) and (\ref{1001}), there is a function $h_1\in\mathcal{H}(\overline{V_1},\alpha)\cap\mathcal{H}(\overline{V_2},\alpha)\cap\mathcal{H}(\overline{V_1*V_2},\alpha)$ such that
\begin{eqnarray*}
\overline{V_1*V_2}(x)\sim \alpha\int_{h_1(x)}^{x-h_1(x)}\overline{V_1}(x-y)\overline{V_2}(y)dy.
\end{eqnarray*}
\textcolor[rgb]{1.00,0.00,0.00}{Then} by (\ref{170616-thm-50}), it is easy to verify that
\begin{eqnarray*}
\overline{V_1*V_2}(x)\sim \alpha^2\overline{G_1*G_2}(x)/a^2(\varepsilon)\sim\alpha^3 e^{-\alpha x}f_{01}\otimes f_{02}(x)/a^2(\varepsilon)=\alpha e^{-\alpha x}f_1\otimes f_2(x).
\end{eqnarray*}

Second, for $i=1,2$, we consider the situation that $V_i$ on $\mathbb{R}$ is lattice with step $1$. Then
\begin{eqnarray*}
\overline{V_i}(k)\sim e^{-\alpha k}l_i(e^k)=e^{-\alpha k}f_i(k)\ \ \text{as}\ k\to\infty,
\end{eqnarray*}
where $f_i\in\mathcal{L}_{0}$ on $\mathbb{N}$.
We linearly extend the function to $\mathbb{R}$ as follows:
$$\textcolor[rgb]{1.00,0.00,0.00}{f_i^{*}(x)}=\sum_{k=1}^\infty \Big(\big(f_i(k+1)-f_i(k)\big)x+(k+1)f_i(k)-kf_i(k+1)\Big)\textbf{1}(k\leq x<k+1),$$
\textcolor[rgb]{1.00,0.00,0.00}{Now we} prove that $e^{-\alpha x}f_i^{*}(x)\downarrow0 $ eventually. \textcolor[rgb]{0.00,0.07,1.00}{Since  $f_i\in\mathcal{L}_{0}$,
$$f_i^{*}(k+1)-f_i^{*}(k)=o\big(\min\{f_i^{*}(k+1), f_i^{*}(k)\}\big).$$
Because $f_i^{*}$ is a linear function, $f_i^{*}\in\mathcal{L}_{0}$ and}
$$f_i^{*}(k+1)-f_i^{*}(k)=o\big( f_i^{*}(x)\textbf{1}(k\leq x<k+1)\big)\ \text{as}\ k\to\infty.$$
Thus, for $k$ large enough and $k\leq x<k+1$,
$$\textcolor[rgb]{1.00,0.00,0.00}{\frac{d}{dx}\big(e^{-\alpha x}f_i^{*}(x)\big)}=e^{-\alpha x}\big(-\alpha f_i^{*}(x)+f_i^{*}(k+1)-f_i^{*}(k)\big)<0.$$
Therefore, there is a continuous distribution $G_i\in\mathcal{L}(\alpha)$ on $\mathbb{R}$ such that
\begin{eqnarray}\label{170616-thm-52}
\overline{G_i}(x)\sim e^{-\alpha x}f_i^{*}(x)
\end{eqnarray}
and
\begin{eqnarray}\label{170616-thm-53}
\overline{V_i}(k)\sim \overline{G_i}(k)\ \ \text{as}\ k\to\infty.
\end{eqnarray}
By the above proof for continuous distributions and (\ref{170616-thm-52}), we have
\begin{eqnarray}\label{170616-thm-54}
\overline{G_1*G_2}(x)\sim \alpha e^{-\alpha x}f_1^{*}\otimes f_2^{*}(x).
\end{eqnarray}
For the lattice case,
\begin{eqnarray}\label{170616-thm-501}
\overline{V_i}(k)\sim\big(\overline{V_i}(k-1)-\overline{V_i}(k)\big)/(e^{\alpha}-1)=P(Y_i=k)/(e^{\alpha}-1).
\end{eqnarray}
\textcolor[rgb]{0.00,0.07,1.00}{When $k\to\infty$, for any $j=j(k)\in\mathbb{R}^+$ such that $j\to\infty$ and $k-j\to\infty$,
\begin{eqnarray}\label{170616-thm-502}
\int_j^{j+1}\overline{G_2}(k-y)\overline{G_1}(y)dy
\sim\overline{G_2}(k-j)\overline{G_1}(j).
\end{eqnarray}}
Then by (\ref{lemma201}) with $c=0$ and the method of proof in (\ref{170616-thm-42}) and (\ref{1001}), we know that
\begin{eqnarray*}
\overline{V_1*V_2}(k)&\sim&P(Y_1+Y_2> k,Y_1>0,Y_2>0)\\
&\sim&\sum_{j=0}^k\overline{V_2}(k-j)V_1(\{j\})\ \ \text{as}\ k\to \infty.
\end{eqnarray*}
Further,  by (\ref{170616-thm-53})-(\ref{170616-thm-502}), we have
\begin{eqnarray*}
\overline{V_1*V_2}(k)&\sim&(e^{\alpha}-1)\sum_{j=0}^k \overline{V_2}(k-j)\overline{V_1}(j)\\
&\sim&(e^{\alpha}-1)\overline{G_1*G_2}(k)/\alpha\\
&\sim& (e^{\alpha}-1) e^{-\alpha k}f_1\otimes f_2(k)\ \ \text{as}\ k\to \infty,
\end{eqnarray*}
that is the first asymptotic formula in (\ref{170618-corol1-1}) holds for $n=2$.

Third, let $X_1$ be a r.v. with a lattice distribution $V_1\in\textcolor[rgb]{0.00,0.07,1.00}{\mathcal{L}_1(\alpha)}$. Then there is a non-lattice distribution $G_1\in\textcolor[rgb]{0.00,0.07,1.00}{\mathcal{L}_1(\alpha)}$ with the corresponding function $f_1$ such that
\begin{eqnarray*}
\overline{V_1}(k)\sim e^{-\alpha k}f_1(k)\sim \overline{G_1}(k)\ \ \text{as}\ k\to\infty.
\end{eqnarray*}
And let $X_2$ be a r.v. with a non-lattice distribution $V_2\in\textcolor[rgb]{0.00,0.07,1.00}{\mathcal{L}_1(\alpha)}$ and the corresponding function $f_2$.
Using (\ref{lemma201}) and the method of proof in (\ref{170616-thm-42}) and (\ref{1001}), we have
\begin{eqnarray*}
\overline{V_1*V_{2}}(x)&\sim&P(Y_1+Y_{2}> x,Y_1>0,Y_{2}>0)\\
&\sim&\sum_{j=0}^{[x]}\overline{V_{2}}(x-j)V_1(\{j\})\\
&\sim& (e^{\alpha}-1)\sum_{j=0}^{[x]} \overline{V_{2}}(x-j)\overline{G_1}(j)\\
&\sim& (e^{\alpha}-1) e^{-\alpha x}f_1\otimes f_{2}(x).
\end{eqnarray*}

Finally, using induction, we can prove that (\ref{170618-corol1-1}) hold for all $n$. Then by Lemma \ref{lem5}, (\ref{170618-corol1-100}) also follows.

\subsection{\bf Proof of Theorem \ref{corol1}}

To prove Theorem \ref{corol1}, we also need the following two lemmas. The \textcolor[rgb]{0.00,0.07,1.00}{first} lemma is the function analogue of Omey et al. (2017). So we omit its proof.

\begin{Lemma}\label{lem7} Assume that $f_i\in\mathcal{R}_{\gamma_i}$ on $\mathbb{R}^+$ and $\int_0^\infty f_i(y)dy=\infty$ for $i=1,2$.

i) If $\gamma_1,\gamma_2>-1$, then
\begin{eqnarray*}
f_1\otimes f_2(x)\sim xf_1(x)f_2(x)B(\gamma_1+1,\gamma_2+1).
\end{eqnarray*}
Thus $f_1\otimes f_2\in\mathcal{R}_{\gamma_1+\gamma_2+1}$.

ii) If $\gamma_1=\gamma_2=-1$, then
\begin{eqnarray*}
f_1\otimes f_2(x)\sim f_1(x)f_2^I(x)+f_2(x)f_1^I(x)
\end{eqnarray*}
and $f_1\otimes f_2\in\mathcal{R}_{-1}$.

iii)  If $\gamma_1=-1,\gamma_2>-1$, then
\begin{eqnarray*}
f_1\otimes f_2(x)\sim f_2(x)f_1^I(x)
\end{eqnarray*}
and $f_1\otimes f_2\in\mathcal{R}_{\gamma_2}$.
\end{Lemma}

For convenience, we write $g_n=f_1\otimes\cdots\otimes f_n$ for each integer $n\ge1$.

\begin{Lemma}\label{lem8}
For any fixed integer $n\ge1$, if $f_i\in\mathcal{R}_{-1},i=1,2,\cdots,n$, then
\begin{eqnarray}\label{170621-lem8-1}
g_n^I(x)\sim\prod_{i=1}^nf_i^I(x)=\sum_{i=1}^n\int_0^x f_i(y)\prod_{1\le j\neq i\le n}f_j^I(y)dy,
\end{eqnarray}
thus $g_n^I\in\mathcal{R}_0$.
\end{Lemma}
\proof We will prove the conclusion by induction. For $n=2$, since $f_i\in\mathcal{R}_{-1}$,
then $f_i^I\in\mathcal{R}_0$ for $i=1,2$. Hence, by ii) of Lemma \ref{lem7}, we obtain that
\begin{eqnarray*}
&&g_2^I(x)=\int_0^xf_1\otimes f_2(y)dy=\int_0^x\int_0^yf_1(y-t)f_2(t)dtdy\\
&=&\int_0^xf_2(t)\int_0^{x-t}f_1(s)dsdt\\
&=&\int_0^xf_2(t)f_1^I(x-t)dt\big(=f_1^I\otimes f_2(x)\big)\\
&\sim&f_1^I(x)f_2^I(x).
\end{eqnarray*}

On the other hand, it holds that
\begin{eqnarray*}
&&\int_0^xf_1(y)f_2^I(y)dy=\int_0^x f_1(y)\int_0^y f_2(z)dzdy\\
&=&\int_0^x f_2(z)\Big(\int_0^x f_1(y)dy-\int_0^z f_1(y)dy\Big)dz\\
&=&f_1^I(x)f_2^I(x)-\int_0^x f_1^I(z)f_2(z)dy.
\end{eqnarray*}
Thus,
$$f_1^I(x)f_2^I(x)=\int_0^xf_1(y)f_2^I(y)dy+\int_0^x f_1^I(y)f_2(y)dy.$$

Now, we assume (\ref{170621-lem8-1}) holds for $n=k$, so  $g_k^I\in\mathcal{R}_0$.
And we continue to prove that it holds for $n=k+1$.
By the same method as the case $n=2$ and induction hypothesis, we have
\begin{eqnarray*}
g_{k+1}^I(x)&=&\int_0^xg_k\otimes f_{k+1} (y)dy=f_{k+1}\otimes g_k^I(x)\\
&\sim&f_{k+1}^I(x)g_k^I(x)=\prod_{i=1}^{k+1}f_i^I(x).
\end{eqnarray*}

Similarly, we can prove the last equation in (\ref{170621-lem8-1}). Therefore, the lemma is proved.$\hfill\Box$\\

Based on the above two lemmas, we prove Theorem \ref{corol1}.

i) By induction,  (\ref{170618-corol1-2}) follows from Theorem \ref{thm1} and  Lemma \ref{lem7} i).

ii) In order to prove (\ref{170709-corol1-1}), we also use induction method. For $n=2$, by Theorem \ref{thm1} and Lemma \ref{lem7} ii), we get
\begin{eqnarray*}
\overline{V_1*V_2}(x)\sim a e^{-\alpha x}f_1\otimes f_2(x)\sim a e^{-\alpha x}\big( f_1(x)f_2^I(x)+f_1^I(x)f_2(x)\big).
\end{eqnarray*}

Assume (\ref{170709-corol1-1}) holds for $n=k$. Then $g_k=f_1\otimes\cdots\otimes f_k\in\mathcal{R}_{-1}$.  For $n=k+1$, by Lemma \ref{lem7} ii) and Lemma \ref{lem8}, we can get that
\begin{eqnarray*}
g_{k+1}(x)&=&g_k\otimes f_{k+1}(x)\sim g_k(x)f_{k+1}^I(x)+f_{k+1}(x)g_k^I(x)\\
&\sim&\sum_{i=1}^{k+1}f_i(x)\prod_{1\le j\neq i\le n}f_j^I(x).
\end{eqnarray*}
Therefore, (\ref{170709-corol1-1}) holds for $n=k+1$.

iii) By Theorem \ref{thm1} and Lemma \ref{lem7} iii), we have
\begin{eqnarray*}
\overline{V_1*V_2}(x)\sim a e^{-\alpha x}f_2(x)f_1^I(x).
\end{eqnarray*}
So the result follows from induction and Lemma \ref{lem8}.

\textcolor[rgb]{1.00,0.00,0.00}{Finally, the conclusion that $V_1*\cdots *V_{n+1}\in\mathcal{L}_{11}(\alpha)$ is follows from the fact that $f_i$ and $f_i^I$ for $1\le i\le n+1$ are regular varying.}

\section{\bf Applications to risk}
\setcounter{equation}{0}\setcounter{Lemma}{0}\setcounter{thm}{0}\setcounter{remark}{0}

In this section, we first introduce a discrete time risk model with stochastic returns, or with both insurance risk and financial risk. Then, based on the results in the previous section, we respectively give some asymptotic estimates of the aggregate net loss and the ruin probability in the finite time period.

\subsection{\bf Model and results}

For every $i\geq 1$, let r.v. $X_i$ be an insurer's net loss (the total amount of claims less premiums) within time period $(i-1,i]$ with distribution $F_i$ on $\mathbb{R}$. And let r.v. $Y_i$ be the stochastic discount factor (the reciprocal of the stochastic accumulation factor) over the same time period with a distribution $G_i$ on $\mathbb{R}^+$. Usually, for $i\ge1$, $X_i$ and $Y_i$ are known as the insurance risk and financial risk, respectively. Further, assume that r.v.'s $X_i,Y_i,i\ge1$ are mutually independent. Then, at each time $n\ge1$, the stochastic present values of aggregate net losses and their maxima respectively are specified as
\begin{eqnarray*}
S_0=0,~~~S_n=\sum_{i=1}^n X_i\prod_{j=1}^iY_j
\end{eqnarray*}
and
\begin{eqnarray*}
M_n=\max_{0\leq k\leq n}S_k.
\end{eqnarray*}
\textcolor[rgb]{1.00,0.00,0.00}{We} call $P(M_n>x)$ the finite time ruin probability at time $n$ denoted by $\psi(x,n)$, where $x\in\mathbb{R}^+\cup\{0\}$ is the initial capital of the insurer.

For the discrete time model, some basic theoretical issues, especially, the closure under product convolution (that is the distribution of the product of independent r.v.'s) for the subexponeantial class, have been studied, see Cline and Samorodnitsky (1994), Tang (2006), Tang (2008) and Xu et al. \textcolor[rgb]{1.00,0.00,0.00}{(2018)}. Based on these results, some asymptotic formulas for the ruin probability are obtained, see, for example, Tang and Tsitsiashvili (2003, 2004). Further, Chen (2011), Yang and Wang (2013), etc. discuss the ruin probability in the model with some dependent insurance risks and financial risks. Most of the existing works assume that the insurance risk and the financial risk respectively follow the same distribution $F$ and $G$, and that the financial risk is dominated by the insurance risk with a subexponeantial distribution, or more general convolution equivalent distribution, namely $\overline{G}(x)=o\big(\overline{F}(x)\big)$ and $F\in\mathcal{S}(\alpha)$ for some $\alpha\in\mathbb{R}^+\cup\{0\}$. In Li and Tang (2015), the dominating relationship between financial risk and insurance risk is not required, however, the distributions of every convex combination of $F$ and $G$ are required to belong to the class $\mathcal{R}^*(\alpha)\subset\mathcal{R}_2(\alpha)$. In order to remove the restrictions on the dominated relationship and enlarge the range of the corresponding distribution class, Hashorva and Li (2014) give the following result. \\

{\bf Theorem A.} In the above model, if, for every $i\geq1$,
$$\overline{F}_i(x)\sim l_i^*(\ln x)(\ln x)^{\gamma^*-1}x^{-\alpha}\ \text{and}\
\overline{G}_i(x)\sim l_i(\ln x)(\ln x)^{\gamma_i-1}x^{-\alpha}$$
for some positive constants $\alpha,\gamma^*$ and $\gamma_i$, and some slowly
varying functions $l_i^*(\cdot)$ and $l_i(\cdot)$,  then, for every $n\geq 1$, letting $\overline{\gamma}_n=\gamma^*+\sum_{i=1}^n \gamma_i$, we have
\begin{eqnarray*}
\psi(x,n)&\sim& P(S_n>x)\sim P\Big(X_n\prod_{j=1}^nY_j>x\Big)\\
&\sim& \Big(\alpha^n \Gamma(\gamma^*)\prod_{i=1}^n\Gamma(\gamma_i)\Big/\Gamma(\overline{\gamma}_n)\Big)l_n^*(\ln x)(\ln x)^{\overline{\gamma}_n-1}x^{-\alpha}\prod_{i=1}^nl_i(\ln x),
\end{eqnarray*}
\textcolor[rgb]{1.00,0.00,0.00}{where $\Gamma(\cdot)$ is the Gamma-function.}

\begin{remark}\label{remark30}
In the proof of the theorem, for some $n\ge1$ and $1\le i\le n$, one of the key objects of study is
$$\overline{V_i}(x)=\overline{F_i}(e^{x})\sim e^{-\alpha x}x^{\gamma^*-1}l^*_i(x).$$
Clearly, for some $\alpha>0$, $V_i\in\mathcal{L}_{11}(\alpha)$ and $F_i\in\mathcal{R}_{11}(\alpha)$ with the same parameter $\gamma^*-1>-1$. In addition, if $\gamma^*<\gamma_j$ for some $1\le j\le n$, then $\overline{F_i}(x)=o\big(\overline{G_j}(x)\big)$ for all $1\le i\le n$, in other words, the financial risk \textcolor[rgb]{0.00,0.07,1.00}{cannot} be dominated by the insurance risk.
\end{remark}

The above remark shows that Theorem A is of great value. However, there are some interesting issues. For example, can the index $\gamma^*-1$ take different values here for \textcolor[rgb]{0.00,0.07,1.00}{different} insurance risks? And can the $\gamma^*$ take zero? Clearly, in the complex insurance business, it is a more reasonable decision that, the insurance risk and the financial risk follow the different distributions and they all have the bigger scope to be chosen. In this section, based on Theorem \ref{thm1} and Theorem \ref{corol1} of the present paper, we use a method different from Hashova and Li (2014) to positively answer these questions.

As in Theorem A, we assume that all distributions are non-lattice. Thus, this means $a=\alpha$.

\begin{thm}\label{thm3-1}
In the above model, for any $n\ge1$ and all $1\le i\le n$, assume that $F_i$ and $G_i$ belong to the class $\mathcal{R}_{11}(\alpha)$ for some $\alpha>0$. More specifically,
$$\overline{F_i}(x)\sim x^{-\alpha }l_i(x)=x^{-\alpha }f_i(\ln x)\ \text{\textcolor[rgb]{0.00,0.07,1.00}{and}}\ \overline{G_i}(x)\sim x^{-\alpha}
l_i^\star(x)=x^{-\alpha }f_i^\star(\ln x),$$
\textcolor[rgb]{0.00,0.07,1.00}{where functions $l_i$ and $l_i^\star$ belong to the class $\mathcal{R}_0$ and
$$\int_0^\infty f_i(y)dy=\int_0^\infty f_i^\star(y)dy=\infty.$$}
In addition, if
\begin{eqnarray}\label{20170805-thm3-1-0}
f_1(x)=o\big(f_k\otimes f^\star_2\otimes \cdots \otimes f^\star_k(x)\big)\ \text{for all}\ 2\le k\leq n,
\end{eqnarray}
then
\begin{eqnarray}\label{20170805-thm3-1-1}
P(S_n>x)\sim \alpha^n x^{-\alpha}f_n\otimes f^\star_1\otimes \cdots \otimes f^\star_n(\ln x).
\end{eqnarray}
Further, if the condition (\ref{20170805-thm3-1-0}) is replaced by the following stronger term:
\begin{eqnarray}\label{20170805-thm3-1-10}
f_{k-1}(x)=o\big(f_k\otimes f^\star_k(x)\big)\ \text{for all}\ 2\le k\leq n,
\end{eqnarray}
then
\begin{eqnarray}\label{20170805-thm3-1-2}
\psi(x,n)\sim \alpha^n x^{-\alpha}f_n\otimes f^\star_1\otimes \cdots \otimes f^\star_n(\ln x).
\end{eqnarray}
\end{thm}

A more precise result than Theorem \ref{thm3-1} is as follows.

\begin{thm}\label{thm4-1}
Under the conditions of Theorem \ref{thm3-1}, for each $n\ge1$, we have the following further results.

i) If $f_i\in \mathcal{R}_{\gamma_i}$ with $\gamma_i>-1$ and $f^\star_i\in \mathcal{R}_{\gamma^\star_i}$ with $\gamma^\star_i>-1$ for $1\le i\le n$, then
\begin{eqnarray}\label{20170805-thm4-1-1}
&&\psi(x,n)\sim P(S_n>x)\nonumber\\
&\sim&\alpha^nx^{-\alpha}(\ln x)^nf_n(\ln x)\prod_{i=1}^nf^\star_i(\ln x)\prod_{j=1}^nB\Big(\gamma_n+\sum_{k=1}^{j-1}\gamma^\star_k+j,\gamma^\star_j+1\Big).
\end{eqnarray}

ii) If $f_i\in \mathcal{R}_{-1}$ and $f^\star_i\in \mathcal{R}_{-1}$ for $1\le i\le n$, then
\begin{eqnarray*}\label{thm324}
&&\psi(x,n)\sim P(S_n>x)\nonumber\\
&\sim& \alpha^nx^{-\alpha}\Big(f_n(\ln x)\prod_{i=1}^nf_i^{\star I}(\ln x)+f_n^I(\ln x)\sum_{i=1}^{n}f_i^\star(x)\prod_{1\le j\neq i\le n}f_j^{\star I}(\ln x)\Big).
\end{eqnarray*}

iii) If $f_i\in \mathcal{R}_{-1}$ and $f^\star_i\in \mathcal{R}_{\gamma^\star_i}$ with $\gamma^\star_i>-1$ for $1\le i\le n$, then
\begin{eqnarray*}\label{thm322}
&&\psi(x,n)\sim P(S_n>x)\nonumber\\
&\sim&\alpha^nx^{-\alpha}(\ln x)^{n-1}f_n^I(\ln x)\prod_{i=1}^nf^\star_i(\ln x)\prod_{j=1}^{n-1}B\Big(\sum_{k=1}^{j}\gamma^\star_k+j,\gamma^\star_{j+1}+1\Big).
\end{eqnarray*}

iv) If $f_i\in \mathcal{R}_{\gamma_i}$ with $\gamma>-1$ and $f^\star_i\in \mathcal{R}_{-1}$ for $1\le i\le n$, then
\begin{eqnarray*}\label{thm323}
\psi(x,n)\sim P(S_n>x)\sim \alpha^nx^{-\alpha}f_n(\ln x)\prod_{i=1}^nf_i^{\star I}(\ln x).
\end{eqnarray*}
\end{thm}

\begin{remark}\label{remark32}
i) According to the second asymptotic expression in (\ref{170618-corol1-1}), if $\gamma_i=\gamma^*-1$ for $1\le i\le n$, then by Lemma \ref{lem7} i), Lemma \ref{lemma41} and Remark \ref{remark30}, the above two conditions (\ref{20170805-thm3-1-0}) and (\ref{20170805-thm3-1-10}) are automatically satisfied. Therefore, the result of Theorem A is properly included in i) of Theorem \ref{thm4-1}.

ii) Of course, we can consider a more complex case that, $\gamma_i=-1,1\le i\le m_1,\gamma_i>-1,m_1+1\le i\le n_1,\gamma^*_i=-1,n_1+1\le i\le m_2$ and $\gamma_i>-1,m_2+1\le i\le n$, where integers $1\le m_1<n_1<m_2<n$ for some $n\ge4$. We omit its details.
\end{remark}

\subsection{\bf Proofs of the results}

To prove the Theorem \ref{thm3-1}, we need the following lemma.
\begin{Lemma}\label{lemma41}
Assume that $f_i\in\mathcal{L}_0$ on $\mathbb{R^+}$ and $\int_0^\infty f_i(y)dy=\infty$ for $i=1,2,3$. If $f_1(x)=o\big(f_2(x)\big)$, then $f_1\otimes f_3(x)=o\big(f_2\otimes f_3(x)\big)$.
\end{Lemma}
\proof It is not difficult to find that there is a function $h\in\bigcap_{i=1}^3\mathcal{H}(f_i,0)$ such that
\begin{eqnarray*}
f_1\otimes f_3(x)&\sim&\int_{h(x)}^{x-h(x)}f_1(x-y)f_3(y)dy\\
&=&o\Big(\int_{h(x)}^{x-h(x)}f_2(x-y)f_3(y)dy\Big)\\
&=&o\big(f_2\otimes f_3(x)\big).
\end{eqnarray*}
Thus, the lemma is proved.$\hfill\Box$

\textbf{Proof of Theorem \ref{thm3-1}.} We first prove (\ref{20170805-thm3-1-1}) with induction. For $n=1$, it follows from Theorem \ref{thm1} that
\begin{eqnarray*}
P(S_1>x)=P(X_1^+Y_1>x)=P(\ln X_1^++\ln Y_1>\ln x)\sim\alpha x^{-\alpha}f_1\otimes f^\star_1(\ln x).
\end{eqnarray*}
Now we assume by induction that (\ref{20170805-thm3-1-1}) holds for $n=k$. Let $S_{k}^{(2)}=\sum_{i=2}^{k+1}X_i\prod_{j=2}^{k+1}Y_j$, then by induction hypothesis, we have
\begin{eqnarray*}
P(S_{k}^{(2)}>x)\sim \alpha^k x^{-\alpha}f_{k+1}\otimes f^\star_2\otimes \cdots \otimes f^\star_{k+1}(\ln x),
\end{eqnarray*}
which means the distribution of $S_{k}^{(2)}$ belongs to the class $\mathcal{R}(\alpha)$. By (\ref{20170805-thm3-1-0}), it is holds that
$$\overline{F_1}(x)=o\big(P(S_{k}^{(2)}>x)\big).$$
Hence,
\begin{eqnarray*}
P(X_1+S_{k}^{(2)}>x)\sim P(S_{k}^{(2)}>x)\sim\alpha^k x^{-\alpha}f_{k+1}\otimes f^\star_2\otimes \cdots \otimes f^\star_{k+1}(\ln x).
\end{eqnarray*}
Further, by Theorem \ref{thm1} and (\ref{10}), we have
\begin{eqnarray*}
P(S_{k+1}>x)&=&P\big(Y_1(X_1+S_{k}^{(2)})^+>x\big)= P\big(\ln Y_1+\ln (X_1+S_{k}^{(2)})^+>\ln x\big)\\
&\sim&\alpha^{k+1} x^{-\alpha}f_{k+1}\otimes f^\star_1\otimes \cdots \otimes f^\star_{k+1}(\ln x),
\end{eqnarray*}
that is (\ref{20170805-thm3-1-1}) holds for $n=k+1$.

Next, we prove (\ref{20170805-thm3-1-2}). On the one hand, it is obvious that
\begin{eqnarray}\label{20170805-thm3-1-3}
P(M_n>x)\geq P(S_n>x)\sim \alpha^n x^{-\alpha}f_n\otimes f^\star_1\otimes \cdots \otimes f^\star_n(\ln x).
\end{eqnarray}
On the other hand, by (\ref{20170805-thm3-1-10}) and Lemma \ref{lemma41}, we can get that
$$P(S_{i}>x)=o\big(P(S_{i+1}>x)\big)=\cdots=o\big(P(S_{n}>x)\big)$$
for all $1\leq i\le n-1$. Therefore, according to the (\ref{20170805-thm3-1-1}) which has been proved, we have
\begin{eqnarray}\label{20170805-thm3-1-4}
P(M_n>x)\leq \sum_{i=1}^nP(S_i>x)\sim P(S_n>x)\sim \alpha^n x^{-\alpha}f_n\otimes f^\star_1\otimes \cdots \otimes f^\star_n(\ln x).
\end{eqnarray}
Combining (\ref{20170805-thm3-1-3}) and (\ref{20170805-thm3-1-4}), we complete the proof of (\ref{20170805-thm3-1-2}).$\hfill\Box$\\

\textbf{Proof of Theorem \ref{thm4-1}} i) We will prove (\ref{20170805-thm4-1-1}) for $S_n$ by induction. For $n=1$, it follows from Theorem \ref{corol1} that
\begin{eqnarray*}
P(S_1>x)=P(\ln X^++\ln Y>\ln x)\sim\alpha x^{-\alpha}(\ln x)f_1(\ln x) f^\star_1(\ln x)B(\gamma_1+1,\gamma^\star_1+1).
\end{eqnarray*}
Now we assume by induction that (\ref{20170805-thm4-1-1}) holds for $n=k$. Let $S_{k}^{(2)}=\sum_{i=2}^{k+1}X_i\prod_{j=2}^{k+1}Y_j$, then we can have
\begin{eqnarray*}
P(S_{k}^{(2)}>x)\sim \alpha^kx^{-\alpha}(\ln x)^kf_{k+1}(\ln x)\prod_{i=2}^{k+1}f^\star_i(\ln x)\prod_{j=2}^{k+1}B\Big(\gamma_{k+1}+\sum_{s=2}^{j-1}\gamma^\star_s+j-1,\gamma^\star_j+1\Big),
\end{eqnarray*}
which means the distribution of $S_{k}^{(2)}$ belongs to $\mathcal{R}(\alpha)$.
Then by (\ref{20170805-thm3-1-0}), we find that
\begin{eqnarray*}
\frac{\overline{F}_1(x)}{P(S_{k}^{(2)}>x)}&=&\frac{\alpha^{-k}(\ln x)^{-k}f_1(\ln x)}{f_{k+1}(\ln x)\prod_{i=2}^{k+1}f^\star_i(\ln x)\prod_{j=2}^{k+1}B\big(\gamma_{k+1}+\sum_{s=2}^{j-1}\gamma^\star_s+j-1,\gamma^\star_j+1\big)}\\
&\to& 0.
\end{eqnarray*}
 Hence,
\begin{eqnarray*}
&&P\big(X_1+S_{k}^{(2)}>x\big)\sim \alpha^kx^{-\alpha}(\ln x)^kf_{k+1}(\ln x)\\
&&\ \ \ \ \ \ \ \ \ \ \ \ \ \ \ \ \ \ \ \ \ \ \ \ \ \ \ \ \cdot\prod_{i=2}^{k+1}f^\star_i(\ln x)\prod_{j=2}^{k+1}B\Big(\gamma_{k+1}+\sum_{s=2}^{j-1}\gamma^\star_s+j-1,\gamma^\star_j+1\Big).
\end{eqnarray*}
Further, by Theorem \ref{corol1}, the assertion holds for $n=k+1$.

For $M_n$, we only need to prove $P(S_{k}>x)=o\big(P(S_{k+1}>x)\big),1\leq k<n$. By (\ref{20170805-thm3-1-10}), it holds that
\begin{eqnarray*}
\frac{P(S_{k}>x)}{P(S_{k+1}>x)}&=&\frac{f_{k}(\ln x)\prod_{i=1}^{k}f^\star_i(\ln x)\prod_{j=1}^{k}B\big(\gamma_{k}
+\sum_{s=1}^{j-1}\gamma^\star_s+j,\gamma^\star_j+1\big)}{\alpha(\ln x) f_{k+1}(\ln x)\prod_{i=1}^{k+1}f^\star_i(\ln x)\prod_{j=1}^{k+1}B\big(\gamma_{k+1}+\sum_{s=1}^{j-1}\gamma^\star_s+j,\gamma^\star_j+1\big)}\\
&\to& 0.
\end{eqnarray*}
So (i) holds.


Similarly, the rest of the theorem can be obtained. We omit its proof. $\hfill\Box$

\section{\bf Examples}
\setcounter{equation}{0}\setcounter{Lemma}{0}\setcounter{thm}{0}\setcounter{remark}{0}\setcounter{exam}{0}\setcounter{pron}{0}

In this section, we respectively give some specific distributions in the classes $\mathcal{L}_{11}(\alpha)$ and $\mathcal{L}_1(\alpha)\setminus\mathcal{L}_{11}(\alpha)$ for some $\alpha\in\mathbb{R}$ with normal shapes and good properties.

According to Theorem 1.3.1 of Bingham et al. (1987), the function $l\in\mathcal{R}_0$ if and only if it may be written in the form
\begin{eqnarray}\label{4101}
l(x)=c_0(x)\exp\Big\{\int_{a_0}^x\varepsilon_0(y)/ydy\Big\}\ \ \text{for}\ \ x\ge \textcolor[rgb]{1.00,0.00,0.00}{a_0},
\end{eqnarray}
where constant $a_0>1$, function $c_0$ is measurable and $c_0(x)\to c_0\in\mathbb{R}^+$, and $\varepsilon_0(x)\to 0$. Thus,
\begin{eqnarray}\label{4102}
f(x)=l(e^x)=c_0(e^x)\exp\Big\{\int_{a_0}^{e^x}\varepsilon_0(y)/ydy\Big\}=c(x)\exp\Big\{\int_a^x\varepsilon(y)dy\Big\},
\end{eqnarray}
where $a=\ln a_0>0$, $c(x)=c_0(e^x)\to c_0$, and $\varepsilon(x)=\varepsilon_0(e^x)$.

\subsection{\bf On the class $\mathcal{L}_{11}(\alpha)$}

We show that the class $\mathcal{L}_{11}(\alpha)$ for some $\alpha\in\mathbb{R}^+$ is a large class of distributions. For example, \textcolor[rgb]{1.00,0.00,0.00}{in (\ref{1000}),} we may respectively take $l(x)=1,\ (\ln\ln x)^\gamma$ or $(\ln x)^\gamma$, then $f(x)=l(e^x)=1,\ x^\gamma$ or $(\ln x)^\gamma$ for some $\gamma\ge-1$, and correspondingly,
$$\overline{V}(x)\sim e^{-\alpha x}f(x)=e^{-\alpha x},\ e^{-\alpha x}(\ln x)^\gamma\ \text{or}\ e^{-\alpha x}x^\gamma.$$
Clearly, $V\in\mathcal{L}_{11}(\alpha)$.

\subsection{\bf On the class $\mathcal{L}_1(\alpha)\setminus\mathcal{L}_{11}(\alpha)$}

We respectively provide two distributions belonging to the class $\mathcal{L}_1(\alpha)\setminus\mathcal{L}_{11}(\alpha)$ for some $\alpha\in\mathbb{R}^+$ and results corresponding to Theorem \ref{corol1}.
\begin{exam}\label{exam41}
In (\ref{4101}), we take $a_0=e$ and $\varepsilon_0(x)=1/\big(2(\ln x)^{1/2}\big)$. Clearly, $\varepsilon_0(x)\to0$, thus $l\in\mathcal{R}_0$. Further, by (\ref{4102}), we have
$$f(x)=l(e^x)=c_0(e^x)\exp\Big\{\int_1^xy^{-1/2}/2dy\Big\}=c_0(e^x)\exp\{x^{1/2}-1\},$$
$$f(x-t)/f(x)\sim\exp\{(x-t)^{1/2}-x^{1/2}\}\to1\ \ \text{for all}\ t,$$
$$f(xt)/f(x)=l(e^{xt})/l(e^x)\sim\exp\{x^{1/2}(t^{1/2}-1)\}\to\infty\ \ \text{for all}\ t>1,$$
and $\int_0^\infty f(y)dy=\infty$. That is $f\in\mathcal{L}_d$ and $f\notin\mathcal{R}_\gamma$ for all $\gamma\in\mathbb{R}$. Therefore, $V\in\mathcal{L}_1(\alpha)\setminus\mathcal{L}_{11}(\alpha)$ for each $\alpha\in\mathbb{R}^+$ and
$$\overline{V}(x)\sim \exp\{-\alpha x+x^{1/2}-1\}.$$
\end{exam}

\begin{pron}\label{pron41}
\textcolor[rgb]{1.00,0.00,0.00}{Let $V$ be a distribution such that
$$\overline{V}(x)\sim \exp\{-\alpha x+Cx^{\beta}+D\}$$
for some constants $C\in\mathbb{R}^+$, $D\in\mathbb{R}$ and $\beta\in(0,1)$. Then $V\in\mathcal{L}_1(\alpha)\setminus\mathcal{L}_{11}(\alpha)$ for some $\alpha\in\mathbb{R}^+$, and for each} $n\geq 2$, there are functions $\xi_i(\cdot):[0,\infty)\mapsto (0,1)$ for $1\le i\le n-1$ such that
\begin{eqnarray*}
\overline{V^{*n}}(x)\sim \alpha^{n-1}\prod_{k=1}^{n-1}\xi_k^{k-1}(x) x^{n-1}\exp\{-\alpha x+C_{n-1}(x)x^{\beta}+nD\},
\end{eqnarray*}
where $C_{n-1}(x)=C\big(\sum_{i=1}^{n-1}(1-\xi_{i}(x))^\beta\prod_{j=i+1}^{n-1}
\xi_{j}^\beta(x)+\prod_{j=1}^{n-1}\xi_{j}^\beta(x)\big)$.
\end{pron}

\proof  \textcolor[rgb]{1.00,0.00,0.00}{Clearly, $V$ belongs to the class $\mathcal{L}_1(\alpha)\setminus\mathcal{L}_{11}(\alpha)$. Further,} when $n=2$, by Theorem \ref{thm1} and the integral mean value theorem, we have

\begin{eqnarray*}
\overline{V^{*2}}(x)&\sim& \alpha e^{-\alpha x+2D}\int_0^x \exp\{C(x-y)^{\beta}+Cy^{\beta}\}dy\\
&=&\alpha x e^{-\alpha x+2D} \int_0^1 \exp\{Cx^{\beta}\big((1-t)^{\beta}+t^{\beta}\big)\}dt\\
&=&\alpha x\exp\{-\alpha x+C_1(x)x^{\beta}+2D\}
\end{eqnarray*}
where $C_1(x)=C\big((1-\xi_1(x))^\beta+\xi_1^\beta(x)\big)$.  $V\in\mathcal{L}_1(\alpha)\setminus\mathcal{L}_{11}(\alpha)$

Further, we can get the conclusion by induction. $\hfill\Box$

\begin{exam}\label{exam42}
Let $f$ be an function on $\mathbb{R}$ such that
\begin{eqnarray*}
f(x)&=&\sum_{k=0}^\infty\Big(2x\textbf{\emph{1}}(4^k\leq x<2\cdot4^k)+(-3x+10\cdot4^k)\textbf{\emph{1}}(2\cdot4^k\leq x<5\cdot4^k/2)\\
&&+x\textbf{\emph{1}}(5\cdot4^k/2\leq x<3\cdot4^k)+(5x-12\cdot4^k)\textbf{\emph{1}}(3\cdot4^k\leq x<4^{k+1})\Big).
\end{eqnarray*}
Clearly, $e^{-\alpha x}f(x)\downarrow0,\ \int_0^\infty f(y)dy=\infty$,
\begin{eqnarray*}
1=\liminf f(x)/x\le\limsup f(x)/x=2
\end{eqnarray*}
and
\begin{eqnarray*}
l(x)=f(\ln x)&=&\sum_{k=0}^\infty\Big(\ln(2x)\textbf{\emph{1}}(4^k\leq x<2\cdot4^k)+\ln(-3x+10\cdot4^k)\textbf{\emph{1}}(2\cdot4^k\leq x<5\cdot4^k/2)\\
&&+\ln(x)\textbf{\emph{1}}(5\cdot4^k/2\leq x<3\cdot4^k)+\ln(5x-12\cdot4^k)\textbf{\emph{1}}(3\cdot4^k\leq x<4^{k+1})\Big).
\end{eqnarray*}
Thus, $l\in\mathcal{R}_0$ and $f$ is not regularly varying. Therefore, in (\ref{1000}), the corresponding distribution $V$ belongs to the class $\mathcal{L}_1(\alpha)\setminus\mathcal{L}_{11}(\alpha)$ for each $\alpha\in\mathbb{R}^+$.
\end{exam}
\textcolor[rgb]{1.00,0.00,0.00}{For this kind of distribution}, we have the following conclusion without the proof which is similar to (\ref{170618-corol1-2}) of Theorem \ref{corol1}.
\begin{pron}\label{pron42}
Under the conditions of Theorem \ref{thm1}, we further require the distribution $V_i$ belongs to the class $\mathcal{L}_2(\alpha)$ for some $\alpha\in\mathbb{R}^+$ and assume there exist regular-variation functions $f_{0i}$ such that
$$0<c_i=\liminf f_i(x)/f_{0i}(x)\le\limsup f_i(x)/f_{0i}(x)=d_i<\infty$$
for any $n\ge1$ and all $1\le i\le n+1$. Further,
if $f_{0i}\in\mathcal{R}_{\gamma_i}$ and $\gamma_i>-1$ for $1\le i\le n+1$, then
\begin{eqnarray*}
\prod_{i=1}^{n+1}c_i&\le&\liminf\overline{V_1*\cdots *V_{n+1}}(x)
\Big/\prod_{j=1}^{n}B\Big(\sum_{k=1}^j\gamma_k+j,\gamma_{j+1}+1\Big)a^{n} e^{-\alpha x}x^{n}\prod_{i=1}^{n+1}f_{0i}(x)\\
&\le&\limsup\overline{V_1*\cdots *V_{n+1}}(x)
\Big/\prod_{j=1}^{n}B\Big(\sum_{k=1}^j\gamma_k+j,\gamma_{j+1}+1\Big)a^{n} e^{-\alpha x}x^{n}\prod_{i=1}^{n+1}f_{0i}(x)\le\prod_{i=1}^{n+1}d_i.
\end{eqnarray*}
where $a$ is given as Theorem \ref{thm1}.
\end{pron}

\textcolor[rgb]{1.00,0.00,0.00}{According to Proposition \ref{pron41} and Proposition \ref{pron42}, we can also obtain some results corresponding to Theorem \ref{thm3-1} and Theorem \ref{thm4-1} for some distributions in the class $\mathcal{L}_1(\alpha)\setminus\mathcal{L}_{11}(\alpha)$.}

\vspace{1cm}
\noindent{\bf Acknowledgements.} The authors are grateful to two referees and Associate Editor for their valuable comments and suggestions which greatly improve the original version of the paper.


\begin{thebibliography}{99}

\bibitem{AFK2003} Asmussen S., Foss S. and Korshunov D. Asymptotics for sums of
random variables with local subexponential behavior. J. Theor. Probab. 16 (2003), 489-518.

\bibitem{BD1996} Bertoin J. and Doney R. A. Some asymptotic results for
transient random walks. Adv. Appl. Probab., 28 (1996), 207-226.

\bibitem{} Bingham N. H., Goldie C. M. and Teugels, J. L.\emph{ Regular Variation}, Cambridge University
Press, Cambridge, 1987.

\bibitem{} Chen Y. The finite-time ruin probability with dependent insurance and financial risks. J. Appl.
Probab. 48 (2011) 4, 1035-1048.

\bibitem{C1964} Chistyakov V.P. A theorem on sums of independent positive
random variables and its application to branching processes. Theory
Probab. Appl. 9 (1964), 640-648.

\bibitem{CNW1973a} Chover J., Ney P. and Wainger S. Functions of probability
measures. J. Anal. Math. 26 (1973a), 255-302.

\bibitem{CNW1973b} Chover J., Ney P. and Wainger S. Degeneracy properties
of subcritical branching processes. Ann. Probab. 1 (1973b), 663-673.

\bibitem{} Cline D.B.H. and Samorodnitsky G. Subexponentiality of the product of independent random
variables. Stoch. Process. Their Appl. 49 (1994), 75-98.



\bibitem{EG1980} Embrechts P. and Goldie C.M. On closure and factorization properties of subexponential and
related distributions. J. Aust. Math. Soc., Ser. A,  29 (1980) 2, 243-256.

\bibitem{EG1982} Embrechts P. and Goldie C.M. On convolution
tails. Stoch. Process. Appl., 13 (1982), 263-278.

\bibitem{FK2007} Foss S. and Korshunov D. Lower limits and equivalences
for convolution tails. Ann. Probab., 1 (2007), 366-383.

\bibitem{} Hashorva E. and Li J. Ecomor and LCR reinsurance with gamma-like
claims. Insurance: Mathematics and Economics, 53 (2013), 206-215.

\bibitem{} Hashorva E. and Li J. Asymptotics for a discrete-time risk model with the emphasis
on financial risk. Probab. Eng. Inform. Sci. 28 (2014) 4, 573-588.

\bibitem{K1989} Kl\"{u}ppelberg C. Subexponential distributions and
characterization of related classes.  Probab. Theory and Related Fields. 82 (1989), 259-269.

\bibitem{L1989} Leslie J. On the non-closure under convolution of the subexponential family.
J. Appl. Probab., 26 (1989), 58-66.

\bibitem{} Li J. and Tang Q. Interplay of insurance and financial risks in a discrete-time model
with strongly regular variation. Bernoulli, 21 (2015) 3, 1800-1823.

\bibitem{LW2012} Lin J. and Wang Y. New examples of heavy-tailed O-subexponential distributions
and related closure properties. Statist. Probab. Lett., 82 (2012), 427-432.

\bibitem{M1989} Murphree E. S. Some new results on the subexponential class. J. Appl. Probab., 26 (1989), 892-897.

\bibitem{} Omey E., Gulck S.V. and Vesilo R. Semi-heavy tails. Accepted by Lithuanian Mathematical Journal.

\bibitem{} Pakes A. G. Convolution equivalence and infinite divisibility. J. Appl. Probab.,
41 (2004), 407-424.

\bibitem{P1980} Pitman E. J. G. Subexponential distribution functions. J. Austral. Math. Soc., Ser. A, 29 (1980), 337-347.

\bibitem{} Tang Q. From light tails to heavy tails through multiplier. Extremes, 11 (2008), 379-391.

\bibitem{} Tang Q. and Tsitsiashvili G. Precise estimates for the ruin probability
in finite horizon in a discretetime model with heavy-tailed insurance and financial risks. Stoch. Process.
Appl. 108 (2003) 2, 299-325.

\bibitem{} Tang Q. and Tsitsiashvili G. Finite- and infinite-time ruin probabilities in the presence of stochastic
returns on investments. Adv. in Appl. Probab. 36 (2004) 4, 1278-1299.


\bibitem{} Tang, Q. Asymptotic ruin probabilities in finite horizon with subexponential losses and associated
discount factors. Probab. Engrg. Inform. Sci. 20 (2006) 1, 103-113.

\bibitem{WXCY2016} Wang Y., Xu H., Cheng D. and Yu C. The local asymptotic estimation for the supremum
of a random walk. Statistical Papers, 59 (2018), 99-126.

\bibitem{W2008} Watanabe T. Convolution equivalence and distributions
of random sums. Probab. Theory Relat. Fields, 142 (2008), 367-397.

\bibitem{} Yang Y. and Wang Y. Tail behavior of the product of two dependent random variables with applications
to risk theory. Extremes 16 (2013) 1, 55-74.

\bibitem{} Xu H., Cheng F., Wang Y. and Cheng D. A necessary and sufficient condition for the subexponentiality
of product convolution. Advances in Applied Probability, 50 (2018) 1, 57-73.

\end{thebibliography}
\end{document}